\documentclass[hidelinks,onefignum,onetabnum]{siamart220329}

\usepackage{amsfonts,upquote}
\usepackage{graphicx,verbatim}

\headers{Computation of Zolotarev rational functions}{Trefethen and Wilber}

\title{{Computation of Zolotarev rational functions}
\thanks{Submitted to the editors DATE.\ \ The work of the second
author was supported by the National Science Foundation under
Grant No.\ DMS-2410045.}}

\author{Lloyd N. Trefethen\thanks{School of Engineering and Applied Sciences,
Harvard University, Cambridge, MA 02138, USA}
(\email{trefethen@seas.harvard.edu}) \and Heather D. Wilber\thanks{Dept.\ of
Applied Mathematics, University of Washington, Seattle, WA 98195, USA}
(\email{hdw27@uw.edu})}

\def\signEF{\hbox{\rm sign}_{E/F}}
\def\complex{{\mathbb{C}}}
\def\Re{\hbox{\rm \kern .8pt Re\kern .6pt}}
\def\Im{\hbox{\rm \kern .8pt Im\kern .6pt}}

\begin{document}

\maketitle

\begin{abstract}
An algorithm is presented to compute Zolotarev rational functions,
that is, rational functions $r_n^*$ of a given
degree that are as small as possible on one set
$E\subseteq\complex\cup\{\infty\}$ relative to their size on another set
$F\subseteq\complex\cup\{\infty\}$ (the third Zolotarev problem).  Along the way we also approximate the sign
function relative to $E$ and $F$ (the fourth Zolotarev problem).
\end{abstract}



\section{\label{intro}Introduction}
\begin{figure}
\begin{center}
\includegraphics[clip, scale=1.3]{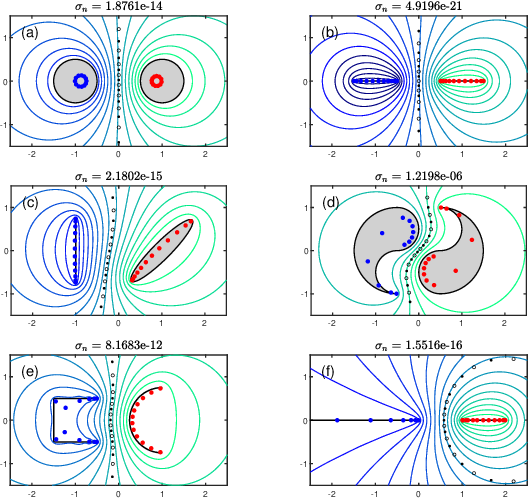}
\end{center}
\caption{\label{fig1}Six examples of degree $n=12$ Zolotarev rational functions
$r_n^*$ defined by connected sets $E$ (on the left in each image) and $F$ (on the right).
The solutions plotted are near-optimal but not exactly so.
The Zolotarev function $r_n^*$ satisfies $|r_n^*(z)| \approx
\min_{t\in F} |r_n^*(t)| = 1$ for $z$ on the boundary of $F$,
and contours show levels $\log_{10} |r_n^*(z)| = -1, -2, \dots$
between the two domains.  Blue dots mark the zeros of\/ $r_n^*$ and red dots
mark the poles.  Black circles and dots mark zeros and poles of the
sign function $\hat r_n^{}$ to be introduced in section~$2$.\ \ The minimum values
$\sigma_n = \|r_n^*\|_E^{}$ are listed in the titles.  Details
of the geometries are given in the appendix.}
\end{figure}

\begin{figure}
\begin{center}
\vskip 10pt
\includegraphics[clip, scale=1.3]{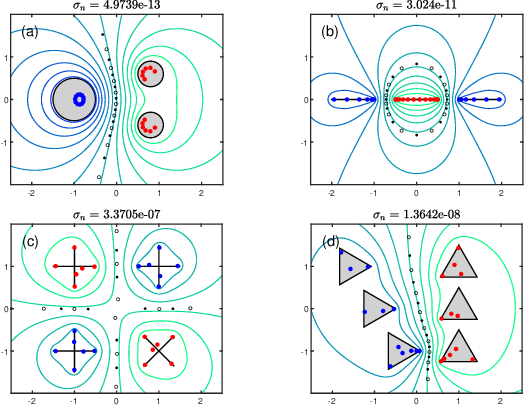}
\end{center}
\caption{\label{fig2}The same as in Figure~$\ref{fig1}$, now
for four problems where $E$ and/or $F$ are disconnected.}
\end{figure}

\begin{figure}
\vskip 10pt
\begin{center}
\includegraphics[clip, scale=1.3]{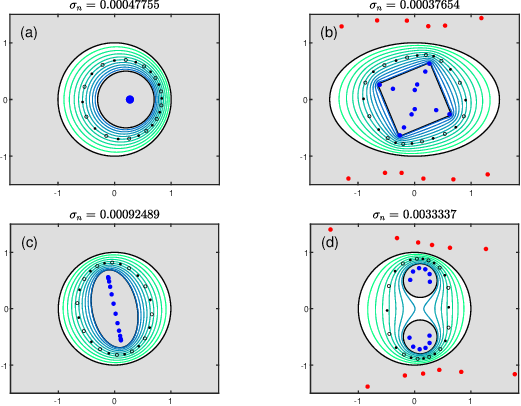}
\end{center}
\caption{\label{fig3}The same as in Figures~$\ref{fig1}$ and\/ $\ref{fig2}$,
now for domains where $F$ encloses $E$.  Since the values of $\sigma_n$ are
not as small in these cases, the contour levels are now at
$\log_{10}|r_n^*(z)| = -1/3, -2/3, \dots.$  Many of the $12$ poles
in each case lie outside these axes.
Mathematically, these
examples are not essentially different from those of
Figures~$\ref{fig1}$ and\/ $\ref{fig2}$, since a 
M\"obius transformation could reduce one type of topology to the other.}
\end{figure}

Figures~\ref{fig1}--\ref{fig3} show 14 examples of computed
Zolotarev rational functions.  The aim of this paper is to present
the mathematics of these functions, explain why they are of interest,
and show how they can be computed numerically by an algorithm
combining AAA-Lawson rational approximation~\cite{aaaL} with the
equivalence of the third and fourth Zolotarev problems of classical
approximation theory~\cite{it}.  Each of the images in the figures
was computed in a fraction of a second on a laptop.  A reliable
method for computing these functions has not been available before.

Let $E$ and $F$ be disjoint closed sets in the complex plane.
For simplicity we assume that $E$ and $F$ each consist of one or a
finite number of continua, such as arcs or domains bounded by arcs.
They should be closed and disjoint in the extended complex plane
$\complex\cup \{\infty\}$, which implies that if $\infty$ belongs
to one of the sets, it does not belong to the other.

For an integer $n\ge 0$, let $R_n$ be the set of rational functions of degree $n$,
which means that any $r\in R_n$ can be written as $p/q$ for some polynomials $p$ and $q$ of
degree at most $n$.
The problem we are concerned with is to find a function $r_n^*\in R_n$
that minimizes the ratio
\begin{equation}
{\max_{z\in E} |r(z)| \over \min_{z\in F} |r(z)|}.
\label{unscaled}
\end{equation}
Since multiplying $r$ by a constant does not change the ratio, we
may normalize the problem by fixing
\begin{equation}
\min_{z\in F} |r(z)| = 1.
\label{norm}
\end{equation}
This gives us what is called the third Zolotarev problem, which might also
be called the {\em Zolotarev ratio problem.}

\medskip
{\em {\bf Problem Z3 = Zolotarev ratio problem.}  Find $r_n^*\in R_n$ with
\rlap{$\min_{z\in F} |r_n^*(z)| = 1$}\hfill\break that attains the minimum
\begin{equation}
\sigma_n = \min_r \|r\|_E^{},
\label{Z3}
\end{equation}
where $\|\cdot\|_E^{}$ denotes the supremum norm over $E$.
}
\medskip

\noindent It is known that a solution $r_n^*$ to (\ref{Z3}) exists, though it need
not be unique, and that $\sigma_n$ satisfies $0 < \sigma_n \le 1$.

Note that there is a certain asymmetry in the normalization
(\ref{norm}), which results in an asymmetry in the contour
lines of Figs.~1--3.  The symmetrical choice would be to scale
$r$ so that the numerator and denominator of (\ref{unscaled})
are reciprocals of each other.  On the other hand (\ref{norm})
is standard, and algebraically simple.

We will say a few words about Figures~\ref{fig1}--\ref{fig3} here,
then turn in section 2 to the mathematical basis of this subject,
the connection between the third and fourth Zolotarev problems of
rational approximation theory.  Section 3 presents our algorithm,
which is based on AAA-Lawson approximation~\cite{aaa} enhanced
by recent modifications for dealing with functions of the flavor
of $f(z) = \hbox{sign}(z)$ and for ensuring convergence of the
Lawson phase.  Sections~4 and~5 discuss applications.
Section 6 finishes with a few closing remarks.

Figures~\ref{fig1}--\ref{fig3}
differ in topology, but all show the same mathematics.  In each case
a close approximation to 
the optimal function $r_n^*$ with $n=12$ has been computed and is
depicted by means of its zeros (blue dots), poles (red dots), and level contours
\begin{equation}
\log_{10} |r_n^*(z)| = -1, -2, \dots ,
\label{levels}
\end{equation}
or in the case of Figure~\ref{fig3},
\begin{equation}
\log_{10} |r_n^*(z)| = -1/3, -2/3, \dots .
\label{levels3}
\end{equation}
The domain $F$ can be recognized as the one containing red poles
and enclosed by green contours.  Because of the normalization
(\ref{norm}), the zero contour $\log_{10} |r_n^*(z)| = 0$ would
enclose $F$ while just touching it; in practice this contour
is very close to the boundary of $F$ (the ``near-circularity
phenomenon''~\cite{starke92}).  The levels in the plots go down  as
far as possible while remaining greater than $\log_{10} \sigma_n$,
implying that they all enclose $E$, which can be recognized as
the region containing blue zeros and enclosed by blue contours.
The colors correspond to the same levels in all the images of
Figures~\ref{fig1} and~\ref{fig2}, and likewise divided by 3 in
the images of Figure~\ref{fig3}.

Each image also shows a chain of black circles and dots, representing
the zeros and poles of the rational function $\hat r_n^{}$ of the
fourth Zolotarev problem, as we will explain in section~2.
Note that about half of our sets $E$ and $F$ have interiors, shaded in
grey in the figures, but our calculations just sample them on
the boundaries.

Looking qualitatively at the figures, we note that smaller values of
$\sigma_n$ correspond to cases where $E$ and $F$ are well separated.
Thus Figure~\ref{fig1}b, with two well separated intervals, gives
$\sigma_n\approx 10^{-20}$, whereas Figure~\ref{fig1}d, where
the yin and yang are interleaved, gives $\sigma_n\approx 10^{-6}$.
In the examples of Figure~\ref{fig3}, $E$ is wholly enclosed by $F$,
and the values of $\sigma_n$ are closer to $1$.

So far as we know, the optimal functions $r_n^*$ in all these
14 example problems are unique (up to multiplication by a
scalar of modulus $1$).  Note that the domains of
Figures~\ref{fig1}d, \ref{fig1}e, \ref{fig1}f, \ref{fig2}a,
\ref{fig2}c, \ref{fig3}a, and \ref{fig3}d each have a line
of symmetry, and those of Figures~\ref{fig1}a, \ref{fig1}b,
\ref{fig2}b, \ref{fig3}b, and \ref{fig3}c each have two lines
of symmetry.  The computed contour lines, poles, and zeros respect
these symmetries quite well, with occasional deviations.  Based on
many experiments with various parameter choices, we believe
that in each image, the value of $\sigma_n$ displayed in the title
is accurate to about two digits or more.

The images convey a vivid impression of electrostatic potential
theory.  It is known that the contour lines can be interpreted as
level curves of the potential generated by positive logarithmic
charges at each pole (i.e., of the form $\log|z-z_k|$) and negative
ones at each zero, and the poles and zeros are positioned in
such a way that the boundaries of $E$ and $F$ are approximate
equipotential surfaces.  For finite $n$, this interpretation
gives bounds on $\sigma_n$, and as $n\to\infty$ it determines the
exponential rate of decrease with $n$.  This theory originates with
Walsh~\cite{walsh}, and for introductions, see~\cite[sec.~6]{ls}
and~\cite[secs.~6--8]{analcont}.

The solutions plotted in Figures~\ref{fig1}--\ref{fig3}, while
nearly optimal, need not be exactly so to plotting accuracy.
In particular, the locations of the poles and zeros may not match
those of the true optimal rational functions (which are only known
analytically in rare cases involving a pair of intervals, a pair
of disks, or a circular annulus~\cite{thesis}).  In Figure~1a, for
example, involving two disks of radius $1/2$ centered at $\pm 1$,
the optimal rational function is known to be a multiple of $[(z-\sqrt
3/2)/(z+\sqrt 3/2)]^{12}$, and the corresponding Zolotarev ratio is
$\sigma_n = [(2-\sqrt 3\kern 1pt)/(2+\sqrt 3\kern 1pt)]^{12} \approx
1.8761\times 10^{-14}$~\cite{starke92}.  This matches the computed
value displayed in the figure, but the computed poles and zeros lie
along small circles rather than coalescing at a point.  Similarly in
Figure~\ref{fig3}a, where $E$ is bounded by the circle of radius
$0.5$ about $0.2$ and $F$ by the unit circle, the exact solution
would have a zero of order 12 at $1/a\approx 0.272$ and a pole of
order 12 at $a \approx 3.68$ with $a = (79/40) + \sqrt{(79/40)^2 -1}$
and $\sigma_n \approx 0.00047755$, again matching the computed value.
Computationally, we find that the zeros are close to the predicted
location but not quite confluent, and the poles, off these axes, lie
approximately on a circle of radius ${\approx\kern 1pt}1$ about the expected location.
To numerical analysts
these are instances of a familiar effect, that small changes in the
values of a polynomial or rational function may be associated with
large changes in its poles and zeros.  To physicists or potential
theorists, they are related to the phenomenon known in potential
theory as {\em balayage,} going back to Poincar\'e or indeed one
might say to Isaac Newton~\cite{gust}.  Outside the unit circle, for
example, a uniform distribution of logarithmic charge on the circle
generates the same potential as a point charge at the origin, and a
finite collection of point charges uniformly spaced along the circle
will approximate the same potential function exponentially closely.

\section{The third and fourth Zolotarev problems}
Yegor Ivanovich Zolotarev (1847--1878) was a student of Chebyshev who
visited Berlin in 1872, where he learned about elliptic functions
from lectures of Weierstrass.  Back in St.\ Petersburg, he applied
these methods to solve a collection of problems involving polynomial
and rational functions posed on two real intervals~\cite{zolo}.
Zolotarev died at age 31 after being hit by a train at the
Tsarskoe Selo station, but his work lived on and was extended,
among others, by Achieser in Kharkiv, Ukraine, who presented
Zolotarev's problems in his approximation theory and elliptic
functions books~\cite{achieser,achieser2}.  Later applications
and generalizations were considered, among others, by Wachspress
at the University of Tennessee~\cite{ew}, Gonchar and colleagues
at the Steklov Institute in Moscow~\cite{gonchar}, Starke at the
University of Karlsruhe~\cite{starke92}, and Istace and Thiran at
the University of Namur~\cite{it}.  Our treatment here follows the
presentation and notation of Istace and Thiran.

We have already stated the third Zolotarev problem Z3, whose
general complex form is due to Gonchar~\cite{gonchar}.  For the fourth Zolotarev
problem, whose complex generalization was
introduced by Istace and Thiran, we first define the {\em sign function} relative to $E$ and $F$:
\begin{equation}
\signEF(z) =
\begin{cases}
-1 & z\in E, \\ +1 & z\in F.
\end{cases}
\end{equation}
(For $z\in\complex\backslash \{E\cup F\}$, $\signEF(z)$ is undefined.)
Problem Z4, which might be called the {\em Zolotarev sign problem,}
is the problem of rational minimax approximation
of $\signEF$ over $E$ and $F$: 

\medskip
{\em {\bf Problem Z4 = Zolotarev sign problem.}  Find $\hat r_n^{}\in R_n$
that attains the minimum
\begin{equation}
\tau_n = \min_r \|r-\signEF\|_{E\cup F}^{},
\label{Z4}
\end{equation}
where $\|\cdot\|_{E\cup F}^{}$ denotes the supremum norm over $E\cup F$.
}
\medskip

\noindent As usual with rational approximation, it is known that
a solution exists, which satisfies $ 0 < \tau_n \le 1$, but it need
not be unique.

In~\cite[p.~143]{achieser2}, Achieser showed that problems Z3 and Z4
are equivalent in the case where $E$ and $F$ are real intervals.
In~\cite{it}, Istace and Thiran showed that the equivalence
generalizes to the complex case.  Here is their theorem, stated
essentially in their words.  The theorem is proved by a direct
calculation, which we do not reproduce.  In particular, it does
not rely on characterizations of solutions to Problems Z3 or Z4.

\begin{figure}
\begin{center}
\includegraphics[clip, scale=0.87]{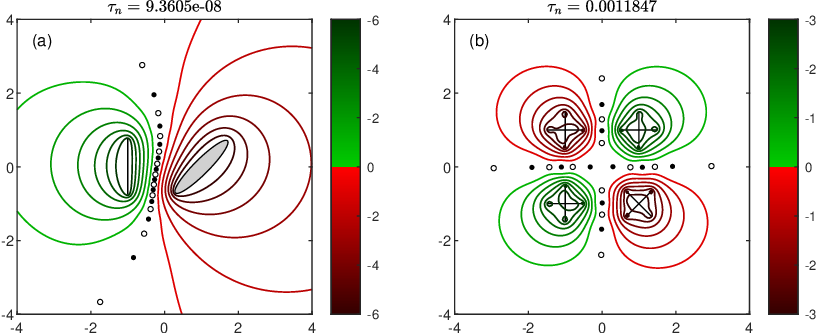}
\end{center}
\caption{\label{fig4} Zeros (\kern .4pt black circles) and poles
(\kern .4pt black dots) of the functions
$\hat r_n^{}$ of Problem Z4 for the examples of Figures~$\ref{fig1}$c and~$\ref{fig2}$c,
together with contour lines showing distances to $+1$ (red) and $-1$ (green).
Specifically, on the left, the red contours show
$\log_{10} |r_n^{}(z)-1| = -1, -2, \dots, -6$ and the green contours show
$\log_{10} |r_n^{}(z)+1| = -1, -2, \dots, -6$;
the pattern on the right is the same except with levels
$-1/2, -1, -3/2,\dots, -3$.
Images like this show that
although Problem Z4 is posed just on the domains $E$ and $F$,
$\hat r_n^{}(z)$ defines an approximation to a sign function throughout the complex plane, with
its poles and zeros delineating an approximate branch cut.
Computation of\/ $\hat r_n^{}$ to solve Problem Z4 is
the first step in our solution of Problem Z3.}
\end{figure}

\begin{theorem}
\label{thm1}
Every solution $r_n^*$ of Problem Z3 is related to a solution 
$\hat r_n^{}$ of Problem Z4 by
\begin{equation}
\hat r_n^{} (z) = {1-\sigma_n\over 1+\sigma_n} \,{r_n^*(z) - \sqrt{\sigma_n}
\over r_n^*(z) + \sqrt{\sigma_n}}, \quad
r_n^*(z) = \sqrt{\sigma_n} \,
{(1-\sigma_n)/(1+\sigma_n) + \hat r_n^{}(z) \over
(1-\sigma_n)/(1+\sigma_n) - \hat r_n^{}(z)}.
\label{relation}
\end{equation}
The minimal values of the two problems satisfy
\begin{equation}
\tau_n = {2\sqrt\sigma_n\over 1 + \sigma_n}, \quad
\sigma_n = \left({\tau_n\over 1 + \sqrt{1-\tau_n^2}\kern .7pt} \right)^2,
\label{minvals}
\end{equation}
and the set of extremal points
\begin{equation}
M = \{z\in E\cup F, ~|\hat r_n^{}(z) - \signEF(z)| = \tau_n\}
\label{extrema}
\end{equation}
is the union of 
$M_1 = \{z\in E, ~|r_n^*(z)| = \sigma_n\}$
and $M_2 = \{z\in F, ~|r_n^*(z)| = 1\}$.
\end{theorem}

Note that in the usual situation $\tau_n, \sigma_n\ll 1$, (\ref{relation}) and (\ref{minvals})
reduce to
\begin{equation}
\hat r_n^{} (z) \approx {r_n^*(z) - \sqrt{\sigma_n}
\over r_n^*(z) + \sqrt{\sigma_n}}, \quad
r_n^*(z) \approx {\tau_n\over 2} \,
{1+ \hat r_n^{}(z) \over 1- \hat r_n^{}(z)}
\label{relation2}
\end{equation}
and
\begin{equation}
\tau_n \approx 2\sqrt\sigma_n, \quad
\sigma_n \approx {\tau_n^2\over 4}.
\label{minvals2}
\end{equation}
The factor $\tau_n/2$ in (\ref{relation2}) is just the scaling
(\ref{norm}), so the
essential point is that rational functions with small ratios
(\ref{unscaled}) come from
approximations $\hat r_n^{}(z)\approx \signEF(z)$
via $(1+\hat r_n^{})/(1-\hat r_n^{})$.

As we will discuss in the next section, our algorithm for solving
Problem Z3 consists of solving Problem Z4 and then transforming
from $\hat r_n^{}$ and $\tau_n$ to $r_n^*$ and $\sigma_n$.  To give
an idea of the mathematics of the equivalence, Figure~\ref{fig4}
shows the poles and zeros of the functions $\hat r_n^{}$ for the
examples of Figures \ref{fig1}c and \ref{fig2}c, which line up along
curves approximating branch cuts for $\signEF$.  The contour lines
show distances to $1$ (red) and $-1$ (green).

\section{Numerical method}
AAA approximation produces a rational function represented in barycentric form.
In the notation of~\cite[eq.~(3.2)]{aaaL}, we have
\begin{equation}
r(z) = \sum_{k=0}^n {\alpha_k\over z-t_k} \left / \sum_{k=0}^n {\beta_k\over z-t_k} \right.,
\label{bary}
\end{equation}
where $t_0,\dots, t_n$ are {\em support points} and
$\alpha_0,\dots,\alpha_n$ and $\beta_0,\dots,\beta_n$ are {\em
barycentric weights}.  Note that $\{t_k\}$ are not poles of
$r$ (assuming the weights are nonzero), but points where the
quotient takes the limiting values $\{\alpha_k/\beta_k\}$.
The zeros and poles of~$r$ are the zeros of the numerator and
denominator of (\ref{bary}), respectively, which can be calculated
accurately by means of a matrix generalized eigenvalue problem
\cite[eq.~(3.11)]{aaa}.

Mathematically, any choice of support points $\{t_k\}$
in (\ref{bary}) would do, but the power of the barycentric
representation lies in its exceptional numerical stability when
the support points are selected in a manner fitted to the function
being represented, as is accomplished by the AAA algorithm.

Following Theorem~\ref{thm1}, we solve Problem Z3 in two steps:
\medskip

(1) {\em Solve Problem Z4 by AAA-Lawson approximation with Chebfun {\tt aaa.m}\kern .7pt ;} 

(2) {\em Convert to a solution of Problem Z3 by $(\ref{relation})$ and $(\ref{minvals})$.}
\medskip

\noindent Step (2) is straightforward, so we just make one comment
on this before turning to discuss the more challenging step (1).
By (\ref{relation}), the zeros and poles of $r_n^*$ are the points
$z$ where $\hat r_n^{} = -p$ and $+p$, respectively, where $p =
(1-\sigma_n)/(1+\sigma_n)$.  To compute these numbers, we take
the barycentric data $\{t_k\}$, $\{\alpha_k\}$, $\{\beta_k\}$
defining $\hat r_n^{}$ and simply subtract or add $p\beta_k$ to each
$\alpha_k$.  Thus $\hat r_n^{}$ is decreased or increased by $p$ at
each support point, hence by the same constant at all $z\in\complex$
since the values at these points determine a degree $n$ rational
interpolant.  To find the zeros and poles of $r_n^*$, it remains only
to perform zerofinding on the barycentric representations of $\hat
r_n^{} \pm p$ in the usual way via generalized eigenvalue problems.

This brings us to the main challenge of our algorithm, step (1)
above, the computation of the degree $n$ rational best approximation
$\hat r_n^{}$ to the sign function $\signEF$ defined by the sets
$E$ and $F$.  The high-level summary is that we find $\hat r_n^{}$
by AAA approximation, which is a fast and robust algorithm for
computing near-best rational approximations~\cite{chebfun,aaa}.  However,
this has proved not as straightforward as one would expect.
Two difficulties arise, and it is because of these that we were
unable to write a paper like the present one a few years ago,
when AAA first became available.  The first difficulty is that AAA
encounters particular challenges when applied to sign functions,
which has not been noticed before.  The second is that for clean
Zolotarev results, it is important to have not just good rational
approximations but nearly optimal ones, requiring the use of the
AAA-Lawson algorithm~\cite{aaaL}, which had a known problem of
non-convergence in certain cases that we have also had to address.

We will discuss these difficulties, and what we have
done about them, in two subsections.

\begin{figure}
\begin{center}
\includegraphics[clip, scale=.9]{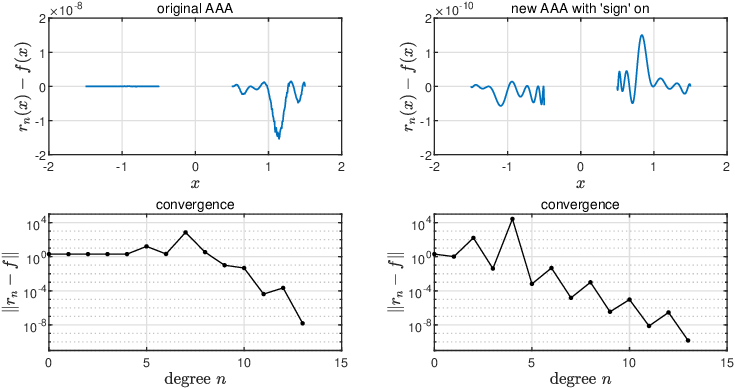}
\end{center}
\caption{\label{fig5}AAA convergence
for degree $n=13$ approximation
$\hat r_n^{} \approx \signEF$ with $E = [-1.5,-0.5]$ and $F = [\kern .5pt 0.5,1.5]$, with each
interval approximated by\/ $200$ Chebyshev points.
On the left, the original AAA iteration reveals difficulties typical of
many $\signEF$ examples.  The ``\kern .5pt blending of singular values''
adjustment described in the text produces the better results on the right.
(This would then further improve to an equioscillatory error curve with the AAA-Lawson
iteration with damping, as illustrated for a more difficult problem in 
Figure~$6$.)
The sawtoothed convergence curve reflects poles falling between sample points at
every other iteration,
a common phenomenon in approximation of even or
odd functions on real domains.}
\end{figure}

\subsection{Modification of AAA for sign functions}
The AAA algorithm has not changed fundamentally since its appearance
in 2018~\cite{aaa}.  Although it has no guarantee of convergence,
we tend to think of it as ``99\% reliable,'' at least when used
with dense enough sample grids and error tolerances well above the
level of noise.  For a wide range of problems, it reliably produces
near-best rational approximations, typically with errors about a
factor of ten above the true minima.

Experiments in approximation of functions of the form $\signEF$,
however, give unsatisfactory results.  The left column of
Figure~\ref{fig5} illustrates what tends to happen.  The initial
iterations achieve nothing, and at the end, the error curve is
very irregular and far from optimal.  After much experimentation we
have discovered the part of the algorithm that causes the trouble.
As described around eq.\ (3.5) of~\cite{aaa}, AAA computes a
singular value decomposition (SVD) at each step to determine a
minimal singular vector defining barycentric weights to solve a
least-squares problem.  The trouble is that in certain cases, the
minimal singular value with standard AAA is degenerate or nearly so,
and this is associated with a singular vector containing zeros or
near-zeros in certain entries.  To fix this, we have introduced a
new option essentially to replace the MATLAB lines \medskip

{\small
\begin{verbatim}
        [~,S,V] = svd(A(J,:),0); 
        wj = V(:,end);
\end{verbatim}
\par}
\medskip

\noindent by the alternative
\medskip

{\small
\begin{verbatim}
        [~,S,V] = svd(A(J,:),0); 
        s = diag(S); wj = V*(1./s.^2); wj = wj/norm(wj);
\end{verbatim}
\par}
\medskip

\noindent
This has the effect that a certain barycentric weight vector {\tt
wj}, rather than being formed from a single minimal singular vector,
is computed as a blend of all the singular vectors with strong bias
towards the minimum.  The effect on most AAA computations appears
to be negligible, but for certain problems such as approximation
of sign functions, there is a real improvement, as suggested in
the right column of Figure~\ref{fig5}.  Note the vertical scales.

The ``strong bias towards the minimum'' in the code
above amounts to an inverse-quadratic weighting of the contribution of
each singular value by the factor $1/s^2$, where $s$ is the corresponding
singular value.  Our initial attempt in this direction used a weighting
by $1/s$, but this led to many approximations being not as close to
minimax as they had been before.  Since most AAA problems don't need
this adjustment at all, it would be unfortunate to introduce a
significant penalty across the board just to improve a small set of
troublesome cases.  After experimentation with various weighting strategies,
we found that $1/s^2$ seemed as good as any, providing great improvement in
sign problems with usually negligible effect in other cases, but we do not claim
that this is precisely optimal in any sense.

In the Chebfun AAA code {\tt aaa.m}, the modification just described was
introduced as an option in July 2024 specified by an optional flag
\verb|'sign'|, which we have invoked for all the computations of
this paper.

\subsection{Modification of Lawson iteration to enhance robustness}
To improve rational approximations from near-best to best, 
the standard method is AAA-Lawson iteration.  As described in
section 3 of~\cite{aaaL}, this is a nonlinear variant of iteratively reweighted
least-squares which often converges to minimax approximations, recognizable by their
equioscillatory error curves in real cases and nearly-circular error curves in
complex problems.  The centerpiece of the Lawson iteration is an adjustment of the
weights in a least-squares problem (linear) according to the errors at the same
points in the current rational approximation (nonlinear).  This is given 
as eq.\ (3.7) of~\cite{aaaL},
\begin{equation}
w_j^{\rm (new)} = w_j |e_j|,
\label{lawson}
\end{equation}
where $w_j$ denotes the current least-squares weight at sample
point $z_j$ and $e_j$ is the current rational approximation error
$e_j = r(z_j)-f_j$.  (Note that these least-squares weights, which
belong to the hundreds or thousands of sample points of the grid,
have nothing to do with the barycentric weights discussed earlier,
which belong just to the subset of $n+1$ barycentric support points.)
In successful cases, iterating (\ref{lawson}) leads to linear
convergence to a weighted least-squares solution that is equal to
the minimax approximation being sought.  Such convergence has long
been known to be guaranteed for the Lawson iteration applied to
linear problems, but with AAA-Lawson, the iteration is nonlinear
and there is no guarantee.

Unfortunately, unlike ordinary AAA iteration, AAA-Lawson has always
been ``just 90\% reliable.''  It fails rather often, and when it
fails, as illustrated in Figure~6.1 of~\cite{aaaL}, the failure
often takes the form of a period-2 oscillation, showing high errors
on the left of a domain at one Lawson step and then high errors
on the right of the domain at the next step.  We have found that
these troubles show up quite often with approximation of $\signEF$
functions.  Zolotarev problems appear to lie in an exceptionally
problematic regime of AAA approximation.

\begin{figure}
\begin{center}
\includegraphics[clip, scale=.9]{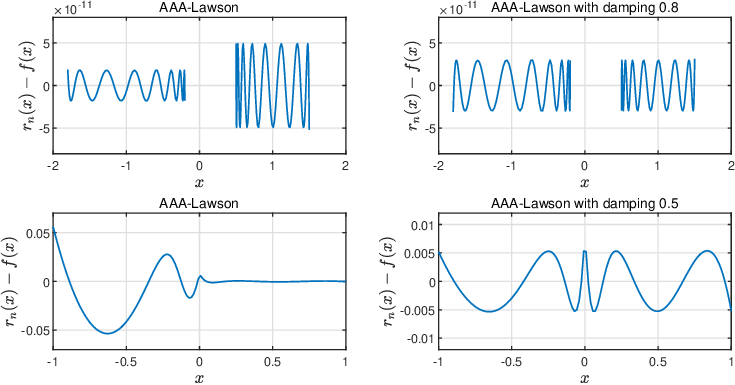}
\end{center}
\caption{\label{fig6}First row: AAA convergence
for degree $15$ approximation of
$\signEF$ with $E = [-1.8,-0.2]$ and $F = [\kern .5pt 0.5,1.5]$, with each
interval approximated by $100$ Chebyshev points.
On the left, the result of AAA followed by $150$ steps of AAA-Lawson iteration.
On the right, the same but with AAA-Lawson applied with damping factor $\delta = 0.8$.
Second row: a similar comparison for degree $n=4$ rational approximation
of\/ $\hbox{\rm ReLU}(x) = \max(x,0)$, sampled in $200$ Chebyshev points, with 
damping factor $\delta = 0.5$.}
\end{figure}

We have investigated this problem and developed a modified
algorithm that often improves matters.  The idea is to replace (\ref{lawson}) by a
modified weight update formula
\begin{equation}
w_j^{\rm (new)} = \left((1-\delta) + {\delta \kern 1pt |e_j|\over \max_j |e_j|}\right) w_j,
\label{lawson2}
\end{equation}
where $\delta\in (0,1]$ is a damping factor.
If $\delta =1$ we have standard AAA-Lawson, whereas
for smaller values we have a more robust iteration that is more likely to converge.
In August 2024 
this modification was
introduced as an option in Chebfun {\tt aaa.m} specified by a
flag \verb|'damping'| followed by the number $\delta$.
Figure~\ref{fig6} gives two illustrations of the improved
convergence of certain iterations with
damping.  For the Zolotarev computations of this paper
we have taken $\delta = 0.95$; details for each
example problem are given in the appendix.

The use of damped iteration as in (\ref{lawson2}) is an old idea
in numerical computation, appearing
in many contexts where oscillations need to be
suppressed.  We are not aware of previous cases of damping in Lawson iterations to minimize
$\infty$-norms as here, which is perhaps understandable since classical Lawson
problems are linear whereas ours is nonlinear.  However, damping
is occasionally applied in $1$-norm applications, where these methods are
more often named iteratively reweighted least-squares (IRLS).
An example is given in~\cite{mold}.

\medskip

In two subsections, we have presented two modifications of the AAA
algorithm that make it more effective in approximating the $\signEF$
functions arising in Zolotarev problems.  We do not regard either of
these adjustments as a definitive solution.  AAA and AAA-Lawson still
give suboptimal results on certain problems,
and we hope that further investigation will
lead to further improvements.  Ideally one would like algorithms
requiring no human intervention such as the specification of a
damping factor $\delta$, and one would like to have a proof that
they always converge.  AAA and AAA-Lawson
continue to advance, but we are a long
way from this state.

Using Chebfun {\tt aaa} and {\tt prz} syntax as of August 2024,
the functions $\hat r_n^{}$ and $r_n^*$ for the example of
Figure~\ref{fig1}a can be computed by this code segment in about
$0.3$ s on our laptop.  The code writes $q$ for $\hat r_n^{}$ and $r$
for $r_n^*$.  \medskip

{\small
\begin{verbatim}
             np = 200; cc = exp(2i*pi*(1:np)'/np);
             E = -1 + .5*cc; F = 1 + .5*cc; 
             fEF = [-ones(size(E)); ones(size(F))];
             [q,qpoles,~,qzeros,zj,fj,wj] = aaa(fEF,[E;F], ...
                  'degree',12,'sign',1,'lawson',200,'damping',0.95);
             tau = norm(fEF-q([E;F]),inf)
             sigma = (tau/(1+sqrt(1-tau^2)))^2
             p = (1-sigma)/(1+sigma);
             r = @(z) sqrt(sigma)*(p+q(z))./(p-q(z));
             [~,~,rpoles] = prz(zj,fj+p,wj);
             [~,~,rzeros] = prz(zj,fj-p,wj);

\end{verbatim}
\par}

\begin{figure}
\begin{center}
\vskip 5pt
\includegraphics[scale=.7]{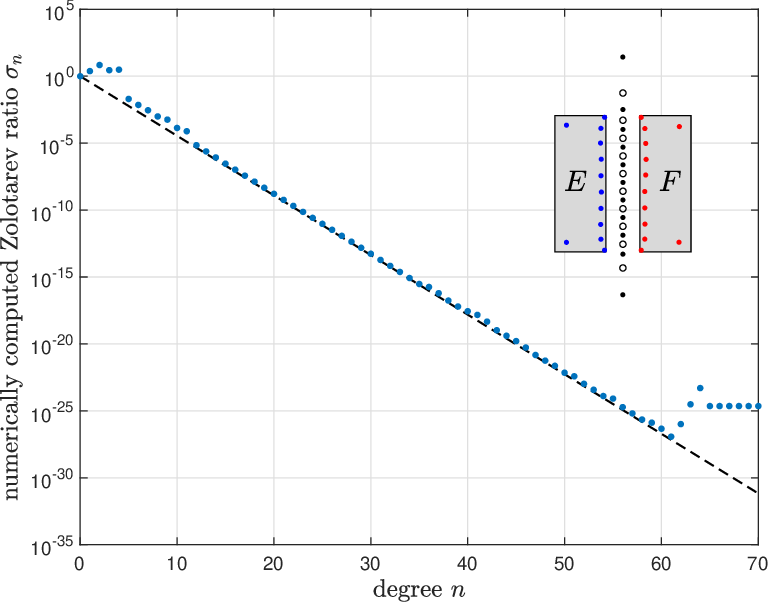}
\end{center}
\caption{\label{rectsfig}Decrease of numerically
computed $\sigma_n$ as a function of $n$ for a problem where
$E$ and $F$ are a pair of rectangles.  
The inset shows the geometry together with
the poles (red) and zeros (blue) of $r_n^*$
and the poles (black dots) and zeros (black circles) of $\hat r_n^{}$
for the case $n=12$.  The dashed line shows the lower bound
$(\ref{capbound})$ from potential theory.
The data points for $n = 1,2,3,4$ are
obviously inaccurate, and
we suspect all the values
for $1\le n\le 11$ should lie closer to the lower bound.}
\end{figure}

To give a further indication of the behavior of our methods, 
Figure~\ref{rectsfig} shows computed $\sigma_n$ as a function 
of $n = 0, 1,\dots, 70$ for a difficult problem involving a pair of rectangles,
as shown in the inset of the figure.  The outer corners
lie at $\pm 1 \pm i$, and the inner corners at $\pm 1/4 \pm i$, and we
took $200$ AAA-Lawson steps with damping parameter $0.95$.
This kind of configuration has been studied in the past, though without
a numerical method available to calculate Zolotarev rational functions
directly~\cite{rtw,thesis}.
Because of the higher degrees involved, this computation took about a minute on our laptop.
The numerical values come close to the lower bound for the Zolotarev
ratio that can be
derived from potential theory via numerical conformal mapping,
\begin{equation}
\sigma_n \ge h^{-n} = e^{-n/\hbox{\scriptsize\rm cap}(E,F)}
\label{capbound}
\end{equation}
where $\hbox{cap}(E,F)$ is what is known as the {\em condenser capacity} of
the pair $E,F$~\cite{ls}.  
We computed the capacity $\hbox{cap}(E,F) \approx 2.78805$ by methods of
rational approximation~\cite{lntcm},
reflecting the curious situation that in the end,
approximate estimates for Zolotarev numbers are not always
simpler to calculate than the Zolotarev numbers themselves.

\section{Applications of the Zolotarev ratio problem Z3}
The most famous application of the Zolotarev ratio problem Z3 is
to ADI-related matrix iterations in numerical linear algebra.
To explain, we begin with the more basic case of polynomial
iterations.

Suppose we want to solve
\begin{equation}
C x = y,
\label{adi1}
\end{equation}
where $C$ is a square matrix, $y$ is a vector, and $x$ is an
unknown vector.  Existence of a unique solution is guaranteed if
$C$ is nonsingular.  Many matrix iterations
generate sequences $x_0, x_1, \dots$ with error vectors
$e_k = C^{-1}y - x_k$ satisfying
\begin{equation}
e_k = p_k(C) e_0,
\label{adi2}
\end{equation}
where $p_k$ is a polynomial of degree $k$ with $p(0) = 1$~\cite{6steps}.
In particular this is true of the Chebyshev, Richardson, 
conjugate gradient, MINRES, and GMRES iterations~\cite{gvl,saad}.
Equation (\ref{adi2}) then implies
\begin{equation}
{\|e_k\|\over \|e_0\|} \le  \|\kern .7pt p_k(C)\|,
\label{adi3}
\end{equation}
where $\|\cdot\|$ is any norm.  If $\|\cdot\|$ is the 2-norm and $\Lambda(C)$ is
the spectrum of $C$, i.e.\ its set of eigenvalues, then (\ref{adi3}) implies
\begin{equation}
{\|e_k\|\over \|e_0\|} \le  \|\kern .7pt p_k\|_\Lambda^{},
\label{adi4}
\end{equation}
where $\|\kern .7pt p_k\|_\Lambda^{} = \max_{z\in\Lambda} |\kern .7pt p_k(z)|$, assuming
$C$ is normal (e.g.\ real symmetric or complex hermitian).
If $C$ is nonnormal, then among various possible generalization of (\ref{adi4}) we have
\begin{equation}
{\|e_k\|\over \|e_0\|} \le  \kappa(V)\|\kern .7pt p_k\|_\Lambda^{},
\label{adi5}
\end{equation}
where $\kappa(V)$ is the 2-norm condition number of any matrix $V$ of
eigenvectors of $C$, assuming one exists.

Equations (\ref{adi1})--(\ref{adi5}) tell us that rapid convergence
of an iteration is guaranteed if it corresponds to a sequence
of polynomials $p_0, p_1,\dots$ normalized by $p_k(0) = 1$ that
converge quickly to $0$ on $\Lambda$.  This connection between
polynomial approximation and matrix iterations has been known and
exploited since the 1950\kern .5pt s.

Rational approximations enter the picture when we generalize (\ref{adi1}) to
the {\em Sylves\-ter equation},
\begin{equation}
AX - XB = Y,
\label{adi6}
\end{equation}
where $A$ is an $m\times m$ matrix, $B$ is an $n\times n$ matrix, $Y$
is an $m\times n$ matrix, and $X$ is an unknown $m\times n$ matrix.
(The special case with $A = -B^*$ and $Y = Y^*$ is called the {\em
Lyapunov equation.}) These equations arise in many applications,
from reduced order modeling and signal processing to the solution
of PDE\kern .5pt s.  See~\cite{sim} and references therein.
Existence of a unique solution is guaranteed if the spectra of
$A$ and $B$ are disjoint.  At first glance (\ref{adi6}) may look
like a very different problem from (\ref{adi1}), and rather niche,
but in fact, (\ref{adi1}) takes the form (\ref{adi6}) in important
cases where the matrix $C$ has special structure.  This observation
became famous with the introduction of {\em Alternating Direction
Implicit} or {\em ADI\/} iteration by Peaceman and Rachford in
1955~\cite{pr}.  In their original model problem, $X$ corresponds
to the unknown values of a discretized finite difference solution
to a two-dimensional linear PDE on an $m\times n$ grid, and the matrices
$A$ and $B$ are finite difference operators with respect to the
two different directions.  Peaceman and Rachford discovered that
(\ref{adi6}) could be solved very efficiently by an iteration in
which one applied $(A-\beta_j I_n)^{-1}$ on the left and $(B-\alpha_j
I_m)^{-1}$ on the right at alternate steps, both involving easy
tridiagonal linear solves, for appropriate constants $\alpha_j$
and $\beta_j$, known as {\em shift parameters.} Omitting details,
which can be found for example in~\cite{sim} and~\cite{thesis},
this leads to an iterative sequence $X_0, X_1, \dots$ with error
matrices $E_k = X - X_k$ satisfying \begin{equation} E_k = r_k(A)
E_0 r_k(B)^{-1}, \label{adi7} \end{equation} where $r_k$ is the
degree $k$ rational function \begin{equation} r_k(z) = \prod_{j=1}^k
{z-\alpha_j\over z-\beta_j}.  \label{adi8} \end{equation} In analogy
to (\ref{adi4}), this implies \begin{equation} {\|E_k\|\over \|E_0\|}
\le {\max_{z\in \Lambda_A} |r(z)| \over \min_{z\in \Lambda_B}
|r(z)|} \label{adi9} \end{equation} if $A$ and $B$ are normal,
with appropriate extensions in the nonnormal case.  And here we
recognize the Zolotarev ratio problem of (\ref{unscaled}).  Thus we
see that the numerical method proposed in this paper offers a new
tool for design and analysis of ADI iterations.

Beyond the basics just outlined, there have been a number of further
applications of the Zolotarev ratio problem in numerical linear
algebra.  One area of application is to rational Krylov
and related methods~\cite{beck,dks}.  Another, investigated recently
by Beckermann and Townsend~\cite{bt}, concerns the case
where the solution matrix $X$ of~(\ref{adi6}) is of large dimension
but of low {\em numerical rank\/}, meaning that its singular values
decay rapidly.  Theorem 2.1 of~\cite{bt} shows show that if $A$ and
$B$ are normal and $Y$ is of rank $\nu \ge 1$ in (\ref{adi6}), then
for each $k$, the singular value of $X$ of index $1+\nu k$ of $X$ can
be bounded in terms of the Zolotarev ratio $\sigma_k$ of (\ref{Z3}),
where $E$ and $F$ are any sets containing the spectra of $A$ and $B$.
Thus when $Y$ is of low rank, $X$ is guaranteed to be of small
numerical rank, with singular values decaying at an exponential
rate determined by the Zolotarev problem.  The tools introduced in
the present paper should make it possible to explore this phenomenon in
much more general cases than have been accessible before.

\section{Applications of the Zolotarev sign problem Z4}
As well as being a step toward the solution of Z3, problem Z4 is
also of interest in its own right.  Problems of this kind arise
in many contexts where one wants to separate one part of a system
computationally from another.  For example, the
1970\kern .5pt s introduced the major technology of
digital signal processing.  A ``recursive'' or ``infinite impulse
response'' low-pass, high-pass, or band-pass filter starts from a
rational function that is nearly constant on one part of the real
axis and nearly zero on another, and finding such a
function is essentially a Zolotarev sign problem~\cite{os}.

Generalizations, as usual, come from numerical linear algebra
and its applications in computational science.  For computing
eigenvalues of large matrices, one of the classes of available
methods is {\em divide and conquer\/} algorithms, where one part
of the spectrum is suppressed relative to another in a possibly
recursive fashion~\cite{bdg,banks}.  This leads quickly to
Zolotarev sign problems, typically starting from the case where
$E$ and $F$ are approximations to the left and right complex
half-planes~\cite{nf}.  In the real case, which goes back to
Zolotarev himself, they may be intervals such as $[-a,-\varepsilon]$
and $[\varepsilon,a]$.  ``Spectral slicing'' ideas of this kind
have found wide generalization through algorithms such as FEAST,
with applications for example in electronic structure calculation
in physics~\cite{feast,lly}.  Here an approximate sign function is
used to isolate a region of the complex plane containing eigenvalues
of mathematical or physical interest.  Mathematically, the issue
is the numerical projection of a large space of functions onto an
interesting smaller-dimensional subspace.

Wherever approximate sign functions are in play, so are rational
approximations to branch cuts, which in turn are close to numerical
quadrature formulas.  Our numerical method for Problem Z4 can
be applied to the derivation of new (and old) quadrature formulas, a topic
to be investigated in a later paper.

\section{Discussion}
Before the current contribution, methods for computing Zolotar\-ev
functions were not available, but various types of approximations have
been discussed in the literature.  We mentioned the lower bound
(\ref{capbound}) in connection with Figure~7,
and to derive upper bounds together with approximate Zolotarev
functions, one can use
Faber rational functions~\cite{rtw} or
Walsh--Fejer, Leja, or generalized Leja
points on the boundary~\cite{ds,starke91}.
All these estimates require some work to apply, however, to solve associated
conformal mapping or optimization problems.

The algorithm we have proposed is not yet bulletproof.  Further
improvements in Zolotarev computations will probably be associated
with further improvements in the AAA and AAA-Lawson algorithms, which we hope
will be stimulated by the discussion here especially in section 3.
For the moment, we have relied on the Chebfun implementation of AAA in
MATLAB/Octave, which
includes the \verb|'sign'| and \verb|'damping'| improvements we have
discussed~\cite{chebfun}.  Implementations of AAA are available in other languages, including
Julia~\cite{driscoll,macmillen} and Python via the SciPy package
(release 1.15.0, January 2025, code written by Jake Bowhay)~\cite{virtanen}.  At present
these lack \verb|'sign'| and \verb|'damping'|,
so they may give less accurate results for Zolotarev problems than
what we have shown here.  Of course, software changes rapidly,
and the details of available options will surely be different in a few years.
What will remain is that approximation algorithms have now advanced
to the point where Zolotarev rational functions can be computed numerically.

\section*{Acknowledgments}
We have benefited from helpful comments from Bernard Beckermann, Nick Hale,
and Yuji Nakatsukasa.  Nakatsukasa in particular contributed to the algorithmic
improvements in AAA described in section~3.

\section*{Appendix.  Details of computed examples} 

For all the examples of Figures~1--3, we ran AAA with \verb|'sign'| on followed by
AAA-Lawson iteration with damping factor $0.95$; the number of iterations
was $200$ in Fig.~1 and $400$ in Figs.~2 and~3.  In the following descriptions
of the approximation sets $E$ and $F$,
$S$ denotes the set of 200 roots of unity, i.e., 200 equispaced
points on the unit circle.

{\em Specifications for Figure $1$.}
(a) Circles $\pm 1 + 0.5S$.
(b) Intervals $[-1.5,-0.5]$ and $[\kern .5pt 0.5,1.5]$, each discretized
by $200$ Chebyshev points.
(c) On the left, $200$ Chebyshev points in the interval
$[-1-0.75\kern .3pt i,-1+0.75\kern .3pt i]$, and on the right, the ellipse 
$1 + (0.2\Re(S) + i\Im(S))/\sqrt{\kern .5pt i\kern 1pt}$.
(d) Let $T$ be the semicircle consisting of the $101$ points of $S$ in
the right half-plane.  The yin figure is composed of three copies of $T$,
two of them reduced to half-size, and all shifted left by $0.5$, and the
yang is the negative of the yin.
(e)  On the right, the semicircle $1-0.74\kern .5pt T$, and on the left, three
sides of a square of side length $0.75$ centered at $-1$, each side
discretized by 100 Chebyshev points.
(f) $200$ Chebyshev points in $[1,2]$ together with
$(-\infty,0\kern .3pt]$ discretized by 200 exponentially graded points: Matlab
\verb|1-logspace(0,5,200)|.

{\em Specifications for Figure $2$.}
(a) Circles $-1+0.5S$ and $0.8+ 0.3S \pm 0.6\kern .3pt i$.
(b) The intervals $[-2,-1]$ $[-0.5,0.5]$, and $[1,2]$, each discretized
by 100 Chebyshev points. 
(c) If $X$ is the cross composed of 100 Chebyshev points in $[-0.5,0.5]$ and
100 Chebyshev points in $[-0.5\kern .5pt i,0.5\kern .5pt i]$, the sets are
composed from $X \pm 1 \pm i$, with the lower-right cross rotated by $\pi/4$. 
(d) Six equilateral triangles scaled to circles of radius $0.5$.  On the right,
the triangles are centered at $1$ and $1\pm i$, and on the left, 
they are rotated by $\pi/6$ and positioned with centers at correspondingly transformed
positions $-1$ and $-1\pm i \mp 1/\sqrt 3$.

{\em Specifications for Figure $3$.}
(a) Circles $S$ and $0.2+0.5S$.
(b) Outside, the ellipse obtained by stretching $S$ by a factor $1.5$ along
the $x$ axis, and inside, the square of side length $1$, each side
discretized by 100 Chebyshev points, rotated by angle $\pi/8$.
(c) Circle $S$ and the ellipse 
$\exp(0.2i)(0.4\Re(S) + 0.7i\Im(S))$.
(d) Circle $S$ and the two smaller circles $0.1 \pm 0.5\kern .3pt i + 0.3S$.

{\em Specifications for Figure $7$.}  The rectangles are discretized by
$50$ Chebyshev points on the ends and $100$ points on the sides.

\end{document}